\documentclass[11pt, reqno]{amsart}
\usepackage{amsmath, amsthm, a4, latexsym, amssymb}

\def\@rmrk#1#2{\refstepcounter
    {#1}\@ifnextchar[{\@yrmrk{#1}{#2}}{\@xrmrk{#1}{#2}}}

\makeatletter\@addtoreset{equation}{section}\makeatother

 \newfont{\bfit}{cmbxti10 scaled 2000}
 \newfont{\biggi}{cmr12 scaled 2000}

 \newcommand{\eps}{\varepsilon}

 \newcommand{\R}{\mathbb{R}}
 
 \newcommand{\N}{\mathbb{N}}

 \newcommand{\prob}{\mathbb{P}}

 \newcommand{\me}{\mathbb{E}}
 
 \renewcommand{\P}{\mathbb{P}}
 
 \newcommand{\one}{\1}

 \newcommand{\skric}{{\mathcal C}}

 \newcommand{\skrie}{{\mathcal E}}
 
 \newcommand{\skrig}{{\mathcal G}}
 \newcommand{\skrih}{{\mathcal H}}

 \newcommand{\skril}{{\mathcal L}}
 \newcommand{\skrim}{{\mathcal M}}

 \newcommand{\skrip}{{\mathcal P}}

 \newcommand{\heap}[2]{\genfrac{}{}{0pt}{}{#1}{#2}}
 \newcommand{\sfrac}[2]{\mbox{$\frac{#1}{#2}$}}

\def\1{{\mathchoice {1\mskip-4mu\mathrm l}      
{1\mskip-4mu\mathrm l}
{1\mskip-4.5mu\mathrm l} {1\mskip-5mu\mathrm l}}}

\newcommand{\eq}{\begin{equation}}
\newcommand{\en}{\end{equation}}

\newenvironment{Proof}
{\vskip0.1cm\noindent{\bf Proof. }{\hspace*{0.3cm}}}%
{\nopagebreak {\hspace*{\fill}\rule{2mm}{2mm}}\\ }

{\nopagebreak {\hspace*{\fill}\rule{2mm}{2mm}}\\ }

\renewcommand{\subsection}{\secdef \subsct\sbsect}
\newcommand{\subsct}[2][default]{\refstepcounter{subsection}
\vspace{0.15cm}
{\flushleft\bf \arabic{section}.\arabic{subsection}~\bf #1  }
\nopagebreak\nopagebreak}
\newcommand{\sbsect}[1]{\vspace{0.1cm}\noindent
{\bf #1}\vspace{0.1cm}}

\newtheorem{theorem}{Theorem}[section]
\newtheorem{lemma}[theorem]{Lemma}
\newtheorem{cor}[theorem]{Corollary}

\newtheoremstyle{thm}{1.5ex}{1.5ex}{\itshape\rmfamily}{}
{\bfseries\rmfamily}{}{2ex}{}

\newtheoremstyle{rem}{1.3ex}{1.3ex}{\rmfamily}{}
{\itshape\rmfamily}{}{1.5ex}{}
\theoremstyle{rem}
\newtheorem{remark}{{\slshape\sffamily Remark}}[]

\refstepcounter{subsection}

\def\thebibliography#1{\section*{References}
  \list%
  {\arabic{enumi}.}
    {\settowidth\labelwidth{[#1]}\leftmargin\labelwidth
    \advance\leftmargin\labelsep
    \parsep0pt\itemsep0pt
    \usecounter{enumi}}
    \def\newblock{\hskip .11em plus .33em minus .07em}
    \sloppy                   
    \sfcode`\.=1000\relax}

\begin{document}
\title[Joint LDP for geometric random  graphs]
{\Large Joint large deviation  result  for empirical measures of the coloured random geometric  graphs}

\author[Kwabena Doku-Amponsah ]{}

\maketitle
\thispagestyle{empty}
\vspace{-0.5cm}

\centerline{\sc By Kwabena Doku-Amponsah }

\centerline{\textit{University of Ghana}}
\centerline{\textit{Statistics  Department, Box LG 115, Legon, Ghana}}
\centerline{\textit{Email: kdoku-amponsah@ug.edu.gh}}

\begin{quote}{\small }{\bf Abstract.}
 We prove joint large deviation principle  for the \emph{ empirical pair measure} and  \emph{empirical locality measure}   of  the \emph{near intermediate} coloured
random  geometric graph models on $n$ points picked uniformly in a  $d-$dimensional torus  of  a  unit  circumference.
From  this result we obtain large
deviation principles for the \emph{number of  edges per vertex}, the \emph{degree distribution and the
proportion of isolated vertices } for the \emph{near intermediate}
random  geometric graph models.
\end{quote}\vspace{0.3cm}


\vspace{0.3cm}

\section{Introduction}

  In this article we study the coloured geometric random graph CGRG, where  $n$ points or vertices or nodes are picked uniformly at random in $[0,1]^d,$ colours or spins are assigned independently from a finite alphabet $\Sigma$ and any two points with colours $a_1,a_2\in\Sigma$ distance at most $r_n(a_1,a_2)$ apart are connected. This random graph models, which has the geometric random graph (see Penrose, 2003) as special case, has been suggested by see (Canning \& Penman, 2003) as a possible extension to the coloured random graph studied in (Biggins \& Penman, 2009), (Doku-Amponsah \& Moerters, 2010), (Doku-Amponsah, 2006),(Bordenave \& Caputo, 2013),(Mukherjee, 2013) and (Doku-Amponsah, 2014).

 The connectivity radius $r_n$ plays similar role
as the connection probability $p_n$ in the  coloured  random
graph  model. Several large  deviation  results  about  the
coloured  random  graphs and  hence Erd\H{o}s-R\'{e}nyi graph have  been  established  recently.
See (O'Connell,1998), (Biggins \& Penman, 2009), (Doku-Amponsah \&  Moerters, 2010), ( Doku-Amponsah, 2006), (Bordenave \& Caputo,2013), (Mukherjee, 2013)
and (Doku-Amponsah, 2014).

Until  recently few  or  no  large  deviation  result  about the CGRG have  been  found. Doku-Amponsah
(2015) proved  joint large  deviation principle  for empirical pair measure  and  the  empirical locality  measure of the CGRG, where $n$
points are uniformly chosen  in $[0,1]^d$, colours or spins are  assigned by drawing without  replacement from  the  pool of, say,  $n\nu_n(a_1)$ colours, and  $n\omega_n(a_1,a_2)$ edges, $a_1,a_2\in\Sigma,$
are randomly inserted  among  the  points for  some  colour  law  $\nu_n:\Sigma\to[0,1]$  and  edge law $\omega_n:\Sigma\times\Sigma\to[0,\infty).$

This  article presents a  full joint large deviation  principle (LDP) for
the empirical pair measure and the  empirical  locality measure  of the CGRG. Refer to
(Doku-Amponsah and Moerters) for  similar result  for  the coloured  random graphs. From  this  large  deviation results  we  obtain LDPs for  graph quantities such as  \emph{number  of  edges per  vertex, the  degree distribution  and the proportion of isolated vertices}  of  geometric random  graphs  in  the  intermediate case.
Our  results  are  similar to those  in  (O'Connell ,1998), (Biggins \& Penman, 2009), (Doku-Amponsah  \&
Moerters, 2010), ( Doku-Amponsah, 2006), (Bordenave  \& Caputo,
  2013), (Mukherjee, 2013)  and (Doku-Amponsah, 2014) for the Erd\"{o}-Renyi graph  except  that the rate  functions  of  the  LDPs  in  our  current setting is  bigger as a  result of  the  effect  of  the  geometric  in  the model.

 As  a  first  step  in  the  proof  of  our  main result, we obtain a joint LDP  for the
\emph{empirical colour measure} and  \emph{empirical pair
measures} for the CGRG, see Theorem~\ref{main},  by  the  exponential change-of-measure techniques and  coupling argument. See  example (Doku-Amponsah \&
 Moerters,  2010) or (Doku-Amponsah,  2006). In  the  next  step, we  use  (Biggins, 2004, Theorem 5(b))  to  mix  Theorem~\ref{main}   and   the  result  (Doku-Amponsah, 2015, Theorem~2.1)  to  obtain  the  full joint LDP for  \emph{empirical  pair  measure  and the  empirical  locality  measure} of  CGRG  model. Refer to  (Doku-Amponsah \&
 Moerters, 2010) or (Doku-Amponsah,  2006) for  further  illustration of  this  method.

 Our  main  motivation  for  studying  this  model  are in two  folds.\\

{\bf Independence testing}: 
     Consider  CGRG  which is  a model  for  Wireless Sensor  Network  as  a very  big  dataset comprising  the  typed sites  and  the  bonds  between  sites. One interesting  question  to  ask  is  how  many  bits  will  be  required  to  code  the  $n$  sites  and  the  bonds  between  sites  with  high probability ? Then,   an  asymptotic equipartition property (AEP) for the  WSN will answer this  question  and  our  LDP  for  the  empirical  measures  of  the  CGRG  will  play a crucial in  the  prove  of  the  AEP.  Further, we can test whether a received codeword $y_n$ of WSN is jointly typical with a candidate sent codeword $x_n$ of  WSN. The probability that two independent sequences $(x_n,y_n)$ ($x_n$ being a codeword other than what was sent when $y_n$ was received) actually appear as dependent is bounded asymptotically as $2^{-nI},$  where the  AEP  is  used  to  obtain  the  bound.   See( Doku-Amponsah,  2016) for  more on  this  application.\\

{\bf Hypothesis testing}:  One of the standard problems in statistics is to decide between two alternative explanations for the data are observed. For example, a transmitter   will send an information on  the  WSN bits by bits in communication systems. There are two possible cases for each transmission: one is that bit 0 of  WSN data is sent (noted as event $H_0$) and the other is that bit 1 of  WSN data  is sent (noted as event $H_1$). In the receiver side, the bit $y$ is to be received as either 0 or 1. Based on the $y$ bit of  WSN data received, we can make a hypothesis whether the event $H_0$ happens (bit 0 was sent at the transmitter) or the event $H_1$ happens (i.e. bit 1 was sent at the transmitter). Of course, we may make mis-judgement, such as we decode that bit 0 was sent but actually bit 1 was sent. We need to make the probability of error in hypothesis testing as low as possible  and the  LDPs  for  CGRG  models  can help  us  specify  the  probability  of  error.\\

In  the remainder of the paper we state and prove  our  LDP results.
In Section~\ref{mainresults}  we   state our LDPs,Theorem~\ref{main2}, Corollary~\ref{Rgg}, Corollary~\ref{RRiv}, Theorem~\ref{main},  and Corollary~\ref{Edgesg}. In Section~\ref{proofmain2} we  present the proof  of  Theorem~\ref{main}.  In Section~\ref{proofmain2} we combine  Theorem~2,1
and (Doku-Amponsah, 2014[b], Theorem~2.1) to obtain the
Theorem~\ref{main2}, using the setup and result of (Biggins, 2004)
to `mix' the LDPs. The paper concludes with the proofs of
Corollary~\ref{ERdd}, Corollary~\ref{RRiv} and  Corollary~\ref{Edgesg} which are given
in Section~\ref{proofcorollaries}.

\section{Statement of the results}\label{mainresults}

\subsection{The joint LDP for empirical pair  measure and  empirical locality  measure  of CGRG.}  In this subsection  we shall look at  a
more general model of random  geometric graphs, the CGCG in
which the connectivity radius depends on  the type or colour or
symbol or  spin of the  nodes. The empirical pair measure and the
empirical locality  measure are our main object of study.

Given   a probability measure $\nu$ on $\Sigma$ and a function
$r_n\colon\Sigma\times\Sigma\rightarrow (0,1]$ we may define the
{\em randomly coloured random geometric graph} or simply
\emph{coloured random geometric graph}~$\skrig$ with $n$ vertices as
follows: Pick vertices  $x_1,...,x_n$  at  random  independently
according to the uniform distribution on $[0,\,1]^d,$ $d\in \N.$ Assign to each
vertex $x_j$ colour $\sigma(x_j)$ independently according to the {\em
colour law} $\nu.$ Given the colours, we join any two vertices
$x_i,x_j$,$(i\not=j)$ by an edge independently of everything else,
if $$\|x_i-x_j\|\le r_n\big[\sigma(x_i),\sigma(x_j)\big].$$ In this  article
we  shall  refer to $r_n(a,b),$   for $a,b\in\Sigma$ as a connection
radius,  and always consider
$$\skrig=(((\sigma(x_i),\sigma(x_j))\,:\,i,j=1,2,3,...,n),E)$$ under the joint law of
graph and colour. We interpret $\skrig$ as  coloured GRG with vertices
$x_1,...,x_n$ chosen at  random uniformly   and independently from
the vertices space $[0,1]^2.$ For  the  purposes of  this  study we
restrict ourselves to  the near intermediate cases .i.e. the
connection radius $r_n$ satisfies the condition $n r_n^d(a,b) \to
C_d(a,b)$ for all $a,b\in \Sigma$, where $C_d\colon\Sigma^2\rightarrow
[0,\infty)$ is a symmetric function, which is not identically equal
to zero.

For any finite or countable set $\Sigma$ we denote by
$\skrip(\Sigma)$ the space of probability measures, and by
$\tilde\skrip(\Sigma)$ the space of finite measures on $\Sigma$,
both endowed with the weak topology.  By  convention we write
$\N=\{0,1,2,...\}.$

We associate with any coloured graph $\skrig$ a probability measure, the
\emph{empirical colour measure}~$\skril^1\in\skrip(\Sigma)$,~by
$$\skril_{\skrig}^{1}(a):=\frac{1}{n}\sum_{j=1}^{n}\delta_{\sigma(x_j)}(a),\quad\mbox{ for $a_1\in\Sigma$, }$$
and a symmetric finite measure, the \emph{empirical pair measure}
$\skril_\skrig^{2}\in\tilde\skrip_*(\Sigma^2),$ by
$$\skril_{\skrig}^{2}(a,b):=\frac{1}{n}\sum_{(i,j)\in E}[\delta_{(\sigma(x_i),\sigma(x_j))}+
\delta_{((\sigma(x_j),\sigma(x_i))}](a,b),\quad\mbox{ for $(a, b)\in\Sigma^2$.
}$$ Note  that the  total mass  
the empirical pair measure is
$2|E|/n$.  Finally we define a further probability measure, the
\emph{empirical neighbourhood measure}
$\skrim_{\skrig}\in\skrip(\Sigma\times\N)$, by
$$\skrim_{\skrig}(a,\ell):=\frac{1}{n}\sum_{j=1}^{n}\delta_{(\sigma(x_i),L(x_j))}(a,\ell),\quad
\mbox{ for $(a,\ell)\in\Sigma\times\N$, }$$ while
$L(x_j)=(l^{x_j}(b),\,b\in\Sigma)$ and $l^{x_j}(b)$ is the number of
vertices  of colour $b$  connected to vertex $x_j$.

For any $\eta\in\skrip(\Sigma\times\N^{\Sigma})$we  denote by $\eta_1$
the $\Sigma-$ marginal of $\eta$ and for  every
$(b,a)\in\Sigma\times\Sigma,$ let $\eta_2$  be the law  of  the pair
$(a,l(b))$ under  the  measure $\eta.$ Define the  measure (finite),
$\langle\eta(\cdot,\ell),\,l(\cdot)\rangle\in\tilde\skrip(\Sigma\times\Sigma)$
by
$$\skrih_2(\eta)(b,a):=
\sum_{l(b)\in\N}\eta_2(a,l(b))l(b), \quad\mbox{ for $a,b\in\Sigma$}$$
and  write $\skrih_1(\eta)=\eta_1.$ We define the function  $\skrih
\colon \skrip(\Sigma\times\N^{\Sigma}) \to \skrip(\Sigma) \times
\tilde\skrip(\Sigma \times \Sigma)$ by
$\skrih(\eta)=(\skrih_1(\eta),\skrih_2(\eta))$ and note that
$\skrih(\skrim_{\skrig})=(\skril_{\skrig}^1, \skril_{\skrig}^2).$ Observe  that
$\skrih_1$ is a continuous function  but $\skrih_2$ is
\emph{discontinuous} in the weak topology. In particular, in  the
summation  $\displaystyle \sum_{l(b)\in\N}\eta_2(a,l(b))l(b)$ the
function $l(b)$ may be unbounded and  so  the  functional
$\displaystyle \eta\to\skrih_2(\eta)$ would not be continuous in the
weak topology. We call a pair of measures
$(\omega, \eta)\in\tilde{\skrip}(\Sigma\times\Sigma)\times\skrip(\Sigma\times\N^{\Sigma})$
\emph{sub-consistent} if
\begin{equation}\label{consistent}
\skrih_2(\eta)(b,a) \le \omega(b,a), \quad\mbox{ for all
$a,b\in\Sigma,$}
\end{equation}
and \emph{consistent} if equality holds in \eqref{consistent}. For a
measure $\omega\in\tilde\skrip_*(\Sigma^2)$ and a measure
$\rho\in\skrip(\Sigma)$, we  recall  from {\bf (Doku-Amponsah \& Moerters, 2010) the rate  function }
$${\mathfrak H}_1(\omega\, \| \, \rho ):=
H\big(\omega\,\|\,C_d\rho\otimes\rho\big)+\|
C_d\rho\otimes\rho \| -\|\omega\|\, ,$$ where the measure
$C_d\rho\otimes\rho\in\tilde\skrip(\Sigma\times\Sigma)$ is defined
by $C_d\rho\otimes\rho(a,b)=C_d(a,b)\rho(a)\rho(b)$ for
$a,b\in\Sigma$. It is not hard to see that $\mathfrak
H_1(\omega\,\|\,\rho)\ge 0$ and equality holds if and only if
$\omega= C_d\rho\otimes\rho$.

For every
$(\omega,\eta)\in\tilde\skrip_*(\Sigma\times\Sigma) \times
\skrip(\Sigma\times\N)$ define a probability measure
$Q_{poi}^{(\omega,\eta)}$ on $\Sigma\times\N$ by
$$Q_{poi}^{(\omega,\eta)}(a\,,\,\ell):=\eta_{1}(a)\prod_{b\in\Sigma} e^{-\frac{\omega(a,b)}{\eta_1(a)}} \,
\frac{1}{\ell(b)!}\,
\Big(\frac{\omega(a,b)}{\eta_1(a)}\Big)^{\ell(b)}, \quad\mbox{for
$a\in\Sigma$, $\ell\in\N$} .$$
We  assume  $d\in\N $ and
  write
  $$\begin{aligned}
\Delta(d)=\left\{
\begin{array}{ll}\sfrac{\pi^{d/2}}{\Gamma\big(\sfrac{(d+2)}{2}\big)}& \mbox{ if  $d\ge 2$}\\
\1& \mbox{ if $d=1,$}
\end{array}\right.
\end{aligned}$$
where  $\Gamma$  is  the  gamma  function.  We now  state the principal theorem
in  this section the LDP for the empirical pair measure and the
empirical locality measure.

\begin{theorem}\label{randomg.LDM}\label{main2}
Suppose that $\skrig$ is a CRGG  with colour law $\nu$ and
connection radii $ r_n\colon\Sigma\times\Sigma\rightarrow[0,1]$
satisfying $n r_n^d(a,b) \to C_d(a,b)$ for some symmetric function
$C\colon\Sigma\times\Sigma\rightarrow [0,\infty)$ not identical to
zero. Then, as $n\rightarrow\infty,$ the pair
$(\skril_{\skrig}^2,\,\skrim_{\skrig})$ satisfies an LDP in
$\tilde{\skrip}_*(\Sigma\times\Sigma)\times\skrip(\Sigma\times\N)$
with good rate function
$$\begin{aligned}
J(\omega,\eta)=\left\{
\begin{array}{ll}H(\eta\,\|\,Q_{poi}^{(\omega,\eta)})+H(\eta_1\,\|\,\nu)+\sfrac{1}{2}{\mathfrak H}_2(\omega\|\eta_1)& \mbox {if $(\omega,\eta)$ consistent  and  $\eta_1=\omega_2,$   }\\
\infty & \mbox{otherwise.}
\end{array}\right.
\end{aligned}$$

$$ {\mathfrak H}_2(\omega\|\eta_1)= \, {\mathfrak H}_1(\omega\,\|\,\eta_1 )-\|\omega\|\log\Delta(d)+(\Delta(d)-\1)\|C_d \eta_1\otimes\eta_1\|.$$
\end{theorem}

\begin{remark}
Note  that  the  first  three  terms  of  the  rate  function is  the  same  as   the  rate  function of  (Doku-Amponsah, 2010,  Theorem~2,1).
  Additionally, the  extra  term  $\sfrac 12( -\|\omega\|\log\Delta(d)+(\Delta(d)-\1)\|C_d \eta_1\otimes\eta_1\|)$    is  positive   and  is  as  a  result  of  the  geometric $[0,\,1]^d$   we  have  incorporated in  the  model.  Moreover, on  typical CGRG
 we  have, $\eta_1=\nu,$
$\omega= \Delta(d) C \,\eta_1 \otimes \eta_1$ and
$$\eta(a,\ell)=\nu(a)\prod_{b\in\Sigma}e^{-\Delta(d) C_d(a,b)\nu(b)}\, \frac{(\Delta(d)C_d(a,b)\nu(b))^{\ell(b)}}{\ell(b)!},\qquad\mbox{for
all $(a,\ell)\in\Sigma\times\N$} .$$ 
 Hence, for some  $\eps$  we  $\P\big\{|\skrim_{\skrig}-\eta\|\ge \eps\big\}\to 0$ as  $n\to\infty.$
\end{remark}  We  write  $$\lambda_1(\delta):=\big(\Delta(d)- \1 \big) \sfrac {c}{2}- \sfrac 12 \, \langle
\delta \rangle\, \log \Delta(d)
\Big)$$

\begin{cor}\label{ERdd}\label{Rgg}
Suppose $D$ is the degree distribution of  the random graph
$\skrig(n,r_n),$  where the connectivity radius $r_n\in (0,1]$
satisfies $n r_n^d \to c \in (0,\infty)$. Then ,as $n\to\infty$, $D$
satisfies an LDP in the space $\skrip(\N \cup \{0\})$ with good rate
function
\begin{equation}\label{randomg.ratedeg}
\begin{aligned}
\lambda_2(\delta)= \left\{ \begin{array}{ll}\Big[ H (d\,\|\,q_{\langle \delta\rangle})+\sfrac 12 \, \langle
\delta \rangle\, \log \big( \sfrac{\langle \delta\rangle}{c}
\big)
 - \sfrac 12 \, \langle \delta \rangle  +  \sfrac {c}{2}\Big] + \lambda_1(\delta),
 & \mbox { if $\langle \delta\rangle< \infty,$ }\\[2mm]
\infty & \mbox { if $\langle \delta\rangle= \infty,$ }
\end{array}\right.
\end{aligned}
\end{equation}
where $q_{k}$ is a poisson distribution with parameter~$k,$ 
and
$\langle \delta\rangle:= \sum_{m=0}^{\infty}m\delta(m)$.
\end{cor}

This   rate function  $\lambda_2$  compares  very  well  with  the  rate  function  of (Doku-Amponsah \& Moerters,  Corollary~2.2, 2010) with  the extra  term  $\lambda_1$  accounting  for  the  the  geometric  effect on the CGRG model.

Next  we give  a  similar result as in (O'Connell, 1998), the
LDP for the  proportion of isolated vertices of the  RGG.
$$\xi_1(y) =\big(\Delta(d)-\1\big) c y(1-y/2)  +
(1-y) \big[ \log\big (\sfrac{\1}{\Delta(d)}\big)-
 \sfrac{(\Delta(d)-\1)c (1-y))}{2} \big]$$

 \begin{cor}\label{ERiv}\label{RRiv}
Suppose $D$ is the degree distribution of  the random graph
$\skrig(n,r_n),$ where the connectivity radius $r_n\in (0,1]$
satisfies $n r_n^d \to c \in (0,\infty)$. Then, as $n\to\infty$, the
proportion of isolated vertices, $D(0)$ satisfies an LDP in $[0,1]$
with good rate function $$\xi_2(y)=y \log y + c y(1-y/2) -
(1-y) \big[ \log\big (\sfrac {c}{a}\big)-
 \sfrac{(a-  c (1-y))^2}{2c(1-y)} \big]+\xi_1(y) \, ,$$
 where $a=a(y)$ is the unique positive solution of $1-e^{-a}=\frac {\Delta(d)c}{a}\, (1-y)$.
 \end{cor}

 From  Corollary~\ref{RRiv} we  deduce  that on a typical random  geometric graphs  the
 number  of isolated  vertices will  grow  like  $n
 e^{-\Delta(d)c}.$ Thus,  as  $n\to\infty,$  the  number  of isolated
 vertices in the geometric random graphs  converges  to $n
 e^{-\Delta(d)c}$  in probability. Again,  the  rate  function $\xi_2$  above  compares  very  well  with  the  result  of  (O'Connell, 1998) with  the  extra  term  $\xi_1$  accounting  for the  influence  of  the geometric  plane $[0, 1]^d$  on  the  model.

  \subsection{The joint LDP  for  the empirical  colour  measure  and  empirical  pair  measure of CGRG}\label{secsub1}
  \begin{theorem}\label{randomg.jointL2L1}\label{main}
Suppose that $\skrig$  is a CGRG   with colour law $\nu$ and
connection radii $ r_n\colon\Sigma^2\rightarrow[0,1]$ satisfying $n
r_n^d(a,b) \to C_d(a,b)$ for some symmetric function
$C_d\colon\Sigma^2\rightarrow [0,\infty)$ not identical to zero. Then,
as  $n\rightarrow\infty,$ the pair $(\skril_\skrig^1,\skril_\skrig^2)$
satisfies an LDP  in
$\skrip(\Sigma)\times\tilde{\skrip}_*(\Sigma^2)$ with good rate
function
\begin{equation}\label{randomg.rateL2L1}
I(\eta_1,\omega)=H(\eta_1\,\|\,\nu)+\sfrac{1}{2}{\mathfrak H}_{2}(\omega\,\|\,\eta_1 ),\,
\end{equation}
where
the measure
$C\eta_1\otimes\eta_1 \in\tilde\skrip_*(\Sigma\times\Sigma)$ is defined
by $C\eta_1\otimes\eta_1 (a,b)=C_d(a,b)\eta_1(a)\eta_1(b)$ for
$a,b\in\Sigma.$

\end{theorem}

Further,  we  state  a  Corollary  of  Theorem~\ref{main} below.

\begin{cor}\label{randomg.L1E}\label{Edgesg}
Suppose that $\skrig$is a CGRG  graph  with colour law $\nu$ and
connection radii $ r_n\colon\Sigma^2\rightarrow[0,1]$ satisfying $n
r_n^d(a,b) \to C_d(a,b)$ for some symmetric function
$C_d\colon\Sigma^2\rightarrow [0,\infty)$ not identical to zero. Then,
as $n\rightarrow\infty,$ the number of edges per vertex $|E|/n$  of
$\skrig$satisfies an LDP in $[0,\infty)$ with good rate function
$$\zeta(x)= x \log x - x+ \inf_{y>0} \big\{ \psi(y) - x \log(y) + y \big\}, $$
where $\psi(y) = \inf H(\eta_1 \, \| \, \nu)$ over all probability
vectors $\eta_1$ with  $\sfrac 12 \Delta(d)\eta_1^{T} C \eta_1 =y$.
\end{cor}

\begin{remark}By taking $C_d(a,b)=c$ one will obtain $\psi(y)=0$ for
$y=\sfrac{\Delta(d) }{2}c$, and $\psi(y)=\infty$ otherwise,  which
establishes  that $|E|/n$ obeys an
LDP in $[0,\infty)$ with  good rate function
$$\zeta(x)= x \log x - x+ \inf_{y>0} \big\{\psi(y) - x \log(\sfrac 12 y) + \sfrac 12 y \big\}, $$
where  $\Delta(d)c =y$.

\end{remark}

\section{Proof of Theorem~\ref{main} }\label{proofmain}
\subsection{Change-of-Measure}

For  any two  points  $U_1$  and  $U_2$  uniformly  and
independently chosen  from  the  space $[0,\,1]^d$  write
$$F(t):=\P\Big\{\|U_1-U_2\|\le t\Big\}.$$

Further, given a function $\tilde{f}\colon\Sigma\rightarrow\R$  and
a symmetric function $\tilde{g}\colon \Sigma^2\rightarrow\R$, we
define the constant $U_{\tilde{f}}$ by
$$U_{\tilde{f}}=\log\sum_{a\in\Sigma}e^{\tilde{f}(a)}\nu(a),$$
and the function $\tilde{h}_n\colon\Sigma^2\rightarrow\R$ by
\begin{equation}\label{hdef}
\tilde{h}_n(a,b)=\log\Big[\big(1-F(r_n(a,b))+ F(r_n(a,b))
e^{\tilde{g}(a,b)}\big)^{-n}\Big],
\end{equation}
for $a,b\in\Sigma.$ We use $\tilde{f}$ and $\tilde{g}$ to define
(for sufficiently large $n$) a new coloured random graph  as
follows:
\begin{itemize}
\item To the $n$ points $x_1,x_2,..,x_n$  picked  independently  and  uniformly in  $[0,1]^d$  we assign colours from
$\Sigma$ independently and identically according to the colour law
$\tilde{\nu}$ defined by
$$ \tilde{\nu}(a)=e^{\tilde{f}(a)-B_{\tilde{f}}}\nu(a).$$
\item Given  any two points $x_u,x_v,$ with $x_u$ carrying colour $a$ and
$x_v$ carrying colour $b$, we connect vertex  $x_u$ to vertex $x_v$ with
probability
$$F(\tilde{r}_{n}(a,b))=\frac{F(r_n(a,b))e^{\tilde{g}(a,b)}}{1-F(r_n(a,b))+F(r_n(a,b))e^{\tilde{g}(a,b)}}.$$
\end{itemize}
We denote the transformed law by $\tilde{\prob}.$  We observe that
$\tilde{\nu}$ is a probability measure and that  $\tilde{\prob}$ is
absolutely continuous with respect to $\prob$ as, for any coloured
graph $\skrig=((\sigma(x_j)\colon j=1,2,3,...,n),E)$,
\begin{align}
\frac{d\tilde\prob}{d\prob}(\skrig) & = \prod_{u\in
V}\sfrac{\tilde{\nu}(\sigma(x_u))}{\nu(\sigma(x_u))}\prod_{(u,v)\in
E}\sfrac{F(\tilde{r}_{n}(\sigma(x_u),\sigma(x_v)))}
{F(r_n(\sigma(x_u),\sigma(x_v)))}\prod_{(u,v)\not\in E}\sfrac{1-F(\tilde{r}_n(\sigma(x_u),\sigma(x_v)))}{1-F(r_n(\sigma(x_u),\sigma(x_v)))}\nonumber\\
& = \prod_{u\in
V}\sfrac{\tilde{\nu}(\sigma(x_u))}{\nu(\sigma(x_u))}\prod_{(u,v)\in
E}\sfrac{F(\tilde{r}_n(\sigma(x_u),\sigma(x_v)))}{F(r_n(\sigma(x_u),\sigma(x_v)))}
\times\sfrac{n-nF(r_n(\sigma(x_u),\sigma(x_v)))}{n-nF(\tilde{r}_n(\sigma(x_u),\sigma(x_v)))}\prod_{(u,v)\in\skrie}\sfrac{n-nF(\tilde{r}_n(\sigma(x_u),\sigma(x_v)))}
{n-nF(r_n(\sigma(x_u),\sigma(x_v)))}\nonumber\\
& = \prod_{u\in V}e^{\tilde{f}(\sigma(x_u))-U_{\tilde{f}}}\prod_{(u,v)\in
E}e^{\tilde{g}(\sigma(x_u),\sigma(x_v))}\prod_{(u,v)\in\skrie}{e^{\frac 1n\, \tilde{h}_n(\sigma(x_u),\sigma(x_v))}}\nonumber\\
& = \exp\big( n\langle \skril_{\skrig}^1,
\tilde{f}-U_{\tilde{f}}\rangle+n\langle \sfrac{1}{2}\skril_{\skrig}^2,
\tilde{g}\rangle+n\langle\sfrac{1}{2}\skril_{\skrig}^1\otimes \skril_{\skrig}^1,
\tilde{h}_n\rangle-\langle \sfrac{1}{2}L_{\Delta}^{1},
\tilde{h}_n\rangle \big) , \label{randomg.Itransform}
\end{align}
where  $$L_{\Delta}^{1}=\sfrac{1}{n}\sum_{u\in
V}\delta_{(\sigma(x_u),\sigma(x_u))}.$$ We write $\langle
g,\omega\rangle:=\sum_{a,b\in\Sigma} g(a,b)\omega(a,b)$ for
$\omega\in\tilde{\skrip}(\Sigma^2)$, and $\langle
f,\rho\rangle:=\sum_{a\in\Sigma} f(a)\rho(a)$ for
$\rho\in\skrip(\Sigma)$, and  note that
$$F(r_n(a,b))=\Delta(d)r_n^d(a,b),\,\mbox{ for all $a,b\in\Sigma^2$}.$$ i.e. the  volume of a
$d$-dimensional (hyper)sphere with radius $r(a,b)$ satisfying
$nr_n^d (a,b)\to C_d(a,b).$

The following lemmas will be useful in the  proofs  of  main
Lemmas.

\begin{lemma}[Euler's lemma]\label{randomg.sequence1}
If $nr_n^d (a,b)\to C_d(a,b)$ for every $a,b\in\Sigma$, then
\begin{equation}\label{randomg.sequence}
\lim_{n\rightarrow\infty}\big[1+\alpha
F(r_n(a,b))\big]^{n}=e^{\alpha\Delta(d) C_d(a,b)},\mbox { for all
$a,b\in\Sigma$ and $\alpha\in\R$. }
\end{equation}
\end{lemma}

\begin{Proof}
Observe that, for any $\eps>0$ and for large $n$  we have
$$\Big[1+\sfrac{\alpha \Delta(d)C_d(a,b)-\eps}{n}\Big]^{n}\le
\Big[1+ \alpha F(r_n(a,b))\Big]^{n}\le \Big[1+\sfrac{\alpha\Delta(d)
C_d(a,b)+\eps}{n}\Big]^{n},$$ by the point-wise convergence. Hence by
the sandwich theorem and Euler's formula we get
(\ref{randomg.sequence}).
\end{Proof}

We  write
$$P^{(n)}(\omega) := \prob\big\{\skril_{\skrig}^1=\omega\big\}.$$

\begin{lemma}\label{randomg.uniexpotightness}\label{randomg.tightness}  The  family of
measures $({P}^n \colon n\in\N)$  is  exponentially tight on
$\skrip(\Sigma)$
\end{lemma}

\begin{Proof}
We use  coupling argument,  see the proof  of  (Doku-Amponsah \& Moerters, 2010, Lemma~5.1)
to show that , for every $\theta>0$,  there exists  $N\in\N$  such
that
$$\limsup_{n\to\infty}\frac{1}{n}\P\big\{|E|>nN\big\}\le-\theta.$$
To  begin,  let $c(d)>\max_{a,b\in\Sigma}C_d(a,b)>0$  and  $nr_n^d(c)\to c(d).$ Using  similar
coupling  arguments  as  in see the proof  of
(Doku-Amponsah \& Moerters, 2010,  Lemma~5.1),  we  can define, for  all sufficiently large
$n,$ a  coloured  random  graph  $\tilde{X}$  with
vertices  $x_1,...,x_n$ chosen  uniformly  from  the  vertices
space $[0, 1]^d,$  colour  law $\eta$ and  connectivity probability $p_n=\P\big\{\|x_i-x_j\|\le r_n(c)\big\}= \Delta(d)r_n^d,$  for  all $i\not=j$  such that  any edge  present  in $\skrig$ is also present  in  $\tilde{X}.$  Let  $|\tilde{E}|$ be  the
number of edges  of  $\tilde{X}.$  Using  the  binomial formula  and
Euler's formula,  we  have  that

\begin{equation}
\begin{aligned}\nonumber
\prob\Big\{ |\tilde{E}|\ge n l\Big\}\le
e^{-nl}\me\big[e^{|\tilde{E}|}\big]&
 =e^{-nl}\sum_{k=0}^{\frac{n(n-1)}{2}}e^{k} \left(\heap{n(n-1)/2}{
 k}\right)\Big(p_n\Big)^{k}\Big(1-p_n\Big)^{n(n-1)/2-k}\\
&= e^{-nl}\Big(
1-p_n+ep_n\Big)^{n(n-1)/2} \le
e^{-nl}e^{nc\Delta(d)(e-1+o(1))},
\end{aligned}
\end{equation}
where we  used $np_n=\Delta(d)nr_n^d\to\Delta(d)c$ in the  last  step.
Now given $\theta>0$ choose $N\in\N$  such that $N>
\theta+\Delta(d)c(e-1)$ and observe that, for sufficiently large
$n,$
\begin{equation}\nonumber
\prob\big\{|E|\ge n N\big\}\le \prob\big\{ |\tilde{E}|\ge
nN\big\}\le e^{-n\theta},
\end{equation}
which implies the statement.
\end{Proof}

\subsection{Proof of the upper bound in Theorem~\ref{main}}\\

We denote by $\skric_1$ the space of functions on $\Sigma$ and  by
$\skric_2$ the space of symmetric functions on $\Sigma^2$, and
define
$$\hat{I}({\eta}_{1},\omega )=\sup_{\heap{f\in\skric_{1}}{g\in\skric_{2}}}
\Big\{\sum_{a\in\Sigma}\big(f(a)-U_{f}\big)\eta_{1}(a) +\sfrac{1}{2}
\sum_{a,b\in\Sigma}g(a,b)\omega(a,b)+ \sfrac{\Delta(d)}{2}
\sum_{a,b\in\Sigma}(1-e^{g(a,b)})C_d(a,b)\eta_{1}(a) \eta_{1}(b) \Big\}$$
for~$({\eta}_{1},\omega )\in\skrip(\Sigma)\times \skrip_*(\Sigma^2)$
\begin{lemma}\label{randomg.L2L1uppbound}
For each closed set
$G\subset\skrip(\Sigma)\times\tilde{\skrip}_*(\Sigma^2),$  we  have
$$\limsup_{n\rightarrow\infty}\sfrac{1}{n}\log\prob\big\{(\skril_{\skrig}^1,\skril_{\skrig}^2)\in
F\big\}\le-\inf_{({\eta}_{1},\omega )\in F}\hat{I}({\eta}_{1},\omega ).$$
\end{lemma}

\begin{Proof}
First let $\tilde{f}\in\skric_1$ and $\tilde{g}\in\skric_2$ be
arbitrary. Define $\tilde{\beta}\colon\Sigma^2\rightarrow\R$ by
$$\tilde{\beta}(a,b)=\Delta(d)(1-e^{\tilde{g}(a,b)})C_d(a,b).$$
Observe that, by Lemma~\ref{randomg.sequence1},
$\tilde{\beta}(a,b)=\lim_{n\rightarrow\infty}\tilde{h}_n(a,b)$ for
all $a,b\in\Sigma$, recalling the definition of $\tilde{h}_n$ from
\eqref{hdef}. Hence, by (\ref{randomg.Itransform}), for sufficiently
large~$n$,
\begin{equation}\nonumber
e^{\max_{a\in\Sigma} | \tilde\beta(a,a)| } \ge\int
e^{\langle\frac{1}{2}L_{\Delta}^{1},\,
\tilde{h}_n\rangle}d\tilde{\prob} =\me\Big\{e^{n\langle \skril_{\skrig}^1,
\tilde{f}-U_{\tilde{f}}\rangle+n\langle\frac{1}{2}\skril_{\skrig}^2,
\tilde{g}\rangle+n\langle\frac{1}{2}\skril_{\skrig}^1\otimes \skril_{\skrig}^1,
\tilde{h}_n\rangle}\Big\},
\end{equation}
where $L_{\Delta}^{1}=\sfrac{1}{n}\sum_{u\in
V}\delta_{(\sigma(x_u),\sigma(x_u))}$ and therefore,
\begin{equation}\label{randomg.estP}
\limsup_{n\rightarrow\infty}\sfrac{1}{n}\log\me\Big\{e^{n\langle
\skril_{\skrig}^1,\tilde{f}-U_{\tilde{f}}\rangle+n\langle
\frac{1}{2}\skril_{\skrig}^2,
\tilde{g}\rangle+n\langle\frac{1}{2}\,\skril_{\skrig}^1\otimes \skril_{\skrig}^1 ,
\tilde{h}_n\rangle}\Big\}\le 0.
\end{equation}

Given $\eps>0$ let
$\hat{I}_{\eps}({\eta}_{1},\omega )=\min\{\hat{I}({\eta}_{1},\omega ),{\eps}^{-1}\}-\eps.$
Suppose  that $({\eta}_{1},\omega )\in G$ and observe that
$\hat{I}({\eta}_{1},\omega )>\hat{I}_{\eps}({\eta}_{1},\omega ).$ We now fix
$\tilde{f}\in\skric_1$ and $\tilde{g}\in\skric_2$ such that
$$\langle\tilde{f}-U_{\tilde{f}},\eta_1\rangle+ \sfrac 12\,\langle
\tilde{g},\omega\rangle+ \sfrac 12 \, \langle
\tilde{\beta},\eta_1\otimes\eta_1 \rangle\ge\hat{I}_{\eps}({\eta}_{1},\omega )
.$$ As $\Sigma$ is finite, there exist open neighbourhoods
$B_{\eta_1}^{1}$ and $B_{\omega}^{2}$ of ${\eta}_{1},\omega $ such that
\begin{equation}\nonumber
\inf_{\heap{\tilde{\eta}_1\in B_{\eta_1}^{1}}{\tilde{\omega}\in
B_{\omega}^{2}}}\big\{\langle \tilde{f}-U_{\tilde{f}} ,
\eta_1\rangle+ \sfrac12 \, \langle
\tilde{g},\tilde\omega\rangle+\sfrac 12 \, \langle
\tilde{\beta},\eta_1\otimes\eta_1\rangle\big\}\ge
\hat{I}_{\eps}({\eta}_{1},\omega )-\eps.
\end{equation}
Using Chebyshev's inequality and (\ref{randomg.estP}) we have that
\begin{equation}\begin{aligned}
\limsup_{n\rightarrow\infty}\sfrac{1}{n}&
\log\prob\big\{(\skril_{\skrig}^1, \skril_{\skrig}^2)\in
B_{\eta_1}^{1}\times B_{\omega}^{2}\big\}\\
&\le\limsup_{n\rightarrow\infty}\sfrac{1}{n}\log
\me\Big\{e^{n\langle
\skril_{\skrig}^1,\tilde{f}-U_{\tilde{f}}\rangle+n\langle
\frac{1}{2}\skril_{\skrig}^2,\tilde{g}\rangle+n\langle\frac{1}{2}\skril_{\skrig}^1\otimes
\skril_{\skrig}^1,\tilde{h}_n\rangle}\Big\}-\hat{I}_{\eps}({\eta}_{1},\omega )+\eps\\
&\le -\hat{I}_{\eps}({\eta}_{1},\omega )+\eps.\label{randomg.LDPballs}
\end{aligned}\end{equation}
Now we use Lemma~\ref{randomg.tightness} with $\theta=\eps^{-1},$ to
choose $N(\eps)\in\N$ such that
\begin{equation}\label{unboundedset}
\limsup_{n\rightarrow\infty}\sfrac{1}{n}\log\prob\Big\{|E|> n
N(\eps) \Big\}\le-\eps^{-1}.
\end{equation}
For this $N(\eps),$ define the set $K_{N(\eps)}$ by
$$K_{N(\eps)}=\Big\{({\eta}_{1},\omega )\in\skrip(\Sigma)\times\tilde{\skrip}_*(\Sigma^2):\|\omega\|\le
2N(\eps)\Big\},$$ and recall that $\|\skril_{\skrig}^2\|=2{|E|}/{n}.$ The
set $K_{N(\eps)}\cap F$ is compact and therefore may be covered by
finitely many sets $B_{\eta_{1,r}}^{1}\times B_{\omega_
{r}}^{2},r=1,\ldots ,m$ with $ (\eta_{1,r},\omega_{r})\in F$ for
$r=1,\ldots ,m.$ Consequently,
\begin{equation}\nonumber
\prob\big\{(\skril_{\skrig}^1, \skril_{\skrig}^2)\in
F\big\}\le\sum_{r=1}^{m}\prob\big\{(\skril_{\skrig}^1, \skril_{\skrig}^2)\in
B_{\eta_{1,r}}^{1}\times
B_{\omega_{r}}^{2}\big\}+\prob\big\{(\skril_{\skrig}^1, \skril_{\skrig}^2)\not\in
K_{N(\eps)}\big\}.
\end{equation}
We may now use (\ref{randomg.LDPballs}) and \eqref{unboundedset} to
obtain, for all sufficiently small $\eps>0$,
\begin{equation}\nonumber
\begin{aligned}
\limsup_{n\rightarrow\infty}\sfrac{1}{n}\log\prob\big\{(\skril_{\skrig}^1,
\skril_{\skrig}^2)\in F\big\}& \le \max_{r=1}^{m} \Big(
\limsup_{n\rightarrow\infty}\sfrac{1}{n}\log\prob\big\{(\skril_{\skrig}^1,
\skril_{\skrig}^2)\in
B_{\eta_{1,r}}^{1}\times B_{\omega_{r}}^{2}\big\} \Big) \vee (-\eps)^{-1} \\
& \le \Big( -\inf_{({\eta}_{1},\omega )\in
G}\hat{I}_{\eps}({\eta}_{1},\omega )+\eps \Big)\vee (-\eps)^{-1}.
\end{aligned}
\end{equation}
Taking $\eps\downarrow 0$ we get the desired statement.
\end{Proof}
Next, we express the rate function in term of relative entropies,
see for example (Dembo \& Zeitouni, 1998, 2.15), and consequently show that it
is a good rate function. Recall the definition of the function $I$
from Theorem~\ref{main}.

\begin{lemma}\label{randomg.Vrate}\ \\ \vspace{-.6cm}
\begin{itemize}
\item[(i)] $\hat{I}({\eta}_{1},\omega )=I({\eta}_{1},\omega ),$ for any
$({\eta}_{1},\omega )\in\skrip(\Sigma)\times\tilde{\skrip}_*(\Sigma^2)$,
\item[(ii)] I is a good rate function and
\item[(iii)] ${\mathfrak{H}_2}(\omega\,\|\, \eta_1)\ge 0$ with equality if and only if
$\omega=\Delta(d)C_d\eta_1\otimes\eta_1 .$
\end{itemize}

\end{lemma}
\begin{Proof}
(i) Suppose that $\omega\not\ll\Delta(d) C_d\eta_1\otimes\eta_1 .$
Then, there exists
 $a_0,b_0 \in\Sigma$ with $C\eta_1\otimes\eta_1 (a_0,b_0)=0$ and
 $\omega(a_0,b_0)>0.$ Define $\hat{g}\colon \Sigma^2\rightarrow \R$  by
$$ \hat{g}(a,b)=\log\big[K(\one_{(a_0,b_0)}(a,b)+\one_{(b_0,a_0)}(a,b))+1\big], \mbox{ for
$a,b\in\Sigma$ and  $K>0.$ }$$ For this choice of $\hat{g}$ and
$f=0$ we have
\begin{align*}\nonumber
&\sum_{a\in\Sigma}\big(f(a)-U_{f}\big)\eta_{1}(a) +\sum_{a,b\in\Sigma}
\sfrac 12
\hat{g}(a,b)\omega(a,b)+\sum_{a,b\in\Sigma}\sfrac {\Delta(d)}{2}(1-e^{\hat{g}(a,b)})C_d(a,b)\eta_{1}(a) \eta_{1}(b) \\
&\ge \sfrac {\Delta(d)}{2} \log(K+1)\omega(a_{0},b_{0}) \to \infty,
\qquad \mbox{ for $K\uparrow\infty.$ }
\end{align*}
Now suppose that $\omega\ll C\eta_1\otimes\eta_1 .$ We have
$$\begin{aligned}
\hat{I}(\eta_1,\omega) & = \sup_{f\in\skric_{1}}
\Big\{\sum_{a\in\Sigma}\Big(f(a)-\log\sum_{a\in\Sigma}e^{{f}(a)}\nu(a)\Big)
\, \eta_{1}(a)  \Big\} \\
&\qquad  +  \sfrac {\Delta(d)}{2}\sum_{a,b\in\Sigma}
C_d(a,b)\eta_1(a)\eta_{1}(b)  +  \sfrac{1}{2} \sup_{g\in\skric_{2}}
\Big\{ \sum_{a,b\in\Sigma}g(a,b)\omega(a,b)-
\Delta(d)\sum_{a,b\in\Sigma} e^{g(a,b)}
C_d(a,b)\eta_1(a)\eta_{1}(b) \Big\}.
\end{aligned}$$
By the variational characterization of relative entropy, the first
term equals $H(\eta_1 \, \| \, \nu)$. By the substitution
$h=\Delta(d)e^{g} \, \frac{C_d \eta_1 \otimes \eta_1}{\omega}$ the
last term equals
$$\begin{aligned}
\sup_{\heap{h\in\skric_{2}}{h \ge 0}} & \sum_{a,b\in\Sigma} \Big[
\log \Big( h(a,b) \frac{\omega(a,b)}{\Delta(d)C_d(a,b)\eta_1(a)
\eta_{1}(b) }\Big)
-h(a,b) \Big]  \, \omega(a,b) \\
& = \sup_{\heap{h\in\skric_{2}}{h \ge 0}} \sum _{a,b\in\Sigma} \big(
\log h(a,b) - h(a,b) \big) \, \omega(a,b) +
\sum _{a,b\in\Sigma} \log \Big( \frac{\omega(a,b)}{\Delta(d)C_d(a,b) \eta_1(a) \eta_1(b)} \Big)\, \omega(a,b) \\
& = - \| \omega\|  + H(\omega \, \| \,\Delta(d) C_d \eta_1\otimes
\eta_1 ),
\end{aligned}$$
where we have used $\sup_{x>0} \log x - x = -1$ in the last step.
This yields that $\hat{I}(\eta_1,\omega)={I}(\eta_1,\omega)$.

(ii) Recall from \eqref{randomg.rateL2L1} and the definition of
${{\mathfrak H}_2}$ that $I(\eta_1, \omega)=H(\omega \,\|\,\nu)+
\sfrac12\,
H\big(\omega\,\|\,\Delta(d)C_d\eta_1\otimes\eta_1 \big)+\sfrac{\Delta(d)}{2}\,
\| C_d\eta_1\otimes\eta_1  \| - \sfrac12\, \|\omega\|$. All summands
are continuous in $\eta_1, \omega$ and thus $I$ is a rate function.
Moreover, for all $\alpha<\infty$, the level sets $\{I \le \alpha\}$
are contained in the bounded set
$\{(\eta_1, \omega)\in\skrip(\Sigma)\times\tilde{\skrip}_*(\Sigma^2)\colon
\,{\mathfrak{H}_2}(\omega\,\|\,\eta_1)\le\alpha\}$ and are
therefore compact. Consequently, $I$ is a good rate function.

(iii)  Consider the nonnegative function $\xi(x)=x\log x-x+1$, for
$x> 0$, $\xi(0)=1$, which has its only root in $x=1$. Note that
\begin{align}\label{randomg.calculusnot}
{\mathfrak{H}_2}(\omega\,\|\,\eta_1)=
\left\{\begin{array}{ll}\int\xi\circ g\,\,d(\Delta(d)C_d\omega\otimes
\omega) &\mbox{ if
$g:=\sfrac{d\omega}{d(\Delta(d)C_d\eta_1\otimes\eta_1 )}\ge 0$  exists, }\\
\infty & \mbox{ otherwise.}
\end{array}\right.
\end{align}
Hence ${\mathfrak{H}_2}(\omega\,\|\,\eta_1)\ge 0$, and if
$\omega=\Delta(d)C_d\eta_1\otimes\eta_1 ,$ then
$\xi(\sfrac{d\omega}{d(\Delta(d)C_d\eta_1\otimes\eta_1 }))=\xi(1)=0$
and so
${\mathfrak{H}_2}(\Delta(d)C_d\eta_1\otimes\eta_1 \,\|\,\omega)=0$.
Conversely, if ${\mathfrak{H}_2}(\omega\,\|\,\omega)=0$, then
$\omega(a,b)>0$ implies $C_d\eta_1\otimes\eta_1 (a,b)>0$, which then
implies $\xi\circ g(a,b)=0$ and further $g(a,b)=1$. Hence
$\omega=\Delta(d)C_d\eta_1\otimes \eta_1,$ which completes the proof
of (iii).
\end{Proof}

\subsection{Proof of the lower bound in Theorem~\ref{main} }\\
We  obtain  the  lower  bound  of Theorem~\ref{main}  from  the
upper bound as  follows:
\begin{lemma}\label{randomg.lowbound1}
For every open set
$O\subset\skrip(\Sigma)\times\tilde{\skrip}_*(\Sigma^2),$  we  have
$$\liminf_{n\rightarrow\infty}\sfrac{1}{n}\log\prob\Big\{(\skril_{\skrig}^1,\skril_{\skrig}^2)\in
O\Big\}\ge -\inf_{(\eta_1, \omega)\in O}I(\eta_1, \omega).$$
\end{lemma}

\begin{Proof}
Suppose $(\eta_1, \omega)\in O,$ with  $\omega\ll
\Delta(d)C_d\eta_1\otimes\eta_1 $. Define
$\tilde{f}_{\omega}\colon\Sigma\rightarrow \R$  by
\begin{equation}\label{randomg.S1}
\begin{aligned}\nonumber
\tilde{f}_{\omega}(a)=\left\{\begin{array}{ll}\log\sfrac{\eta_1(a)}{\nu(a)},
&\mbox{if $\eta_1(a)> 0$,  }\\
0, & \mbox{otherwise.}
\end{array}\right.
\end{aligned}
\end{equation}
and $\tilde{g}_{\omega}\colon\Sigma^2\rightarrow \R$ by
\begin{equation}\label{randomg.S2}
\begin{aligned}\nonumber
\tilde{g}_{\omega}(a,b)=\left\{\begin{array}{ll}\log\sfrac{\omega(a,b)}{\Delta(d)C_d(a,b)\eta_1(a)\eta_1(b)},&\mbox{if $\omega(a,b)>0$, }\\
0, & \mbox{otherwise.}
\end{array}\right.
\end{aligned}
\end{equation}
In addition, we let
$\tilde{\beta}_{\omega}(a,b)=\Delta(d)C_d(a,b)(1-e^{{\tilde
g}_{\omega}(a,b)})$ and note that
$\tilde{\beta}_{\omega}(a,b)=\lim_{n\rightarrow\infty}\tilde{h}_{\omega,
n}(a,b),$ for all $a,b\in\Sigma$  where
\begin{equation}\nonumber
\tilde{h}_{\omega,
n}(a,b)=\log\Big[\big(1-F(r_n(a,b))+F(r_n(a,b))e^{\tilde{g}_{\omega}(a,b)}\big)^{-n}\Big].
\end{equation}
Choose $B_{\eta_1}^{1} ,B_{\omega}^{2}$  open neighbourhoods of
$\eta_1, \omega,$ such that  $B_{\eta_1}^{1}, \times
B_{\omega}^{2}\subset O$ and for all
$(\tilde{\omega},\tilde{\omega})\in B_{\eta_1}^{1}\times
B_{\omega}^{2}$
$$\langle \tilde{f}_{\omega},\eta_1\rangle+\sfrac 12\, \langle\tilde{g}_{\omega},\omega\rangle+
\sfrac 12\,
\langle\tilde{\beta}_{\omega},\eta_1\otimes\eta_1 \rangle-\eps\le
\langle \tilde{f}_{\omega},\tilde{\eta}_1\rangle+ \sfrac 12\, \langle
\tilde{g}_{\omega},\tilde{\omega}\rangle+ \sfrac 12\,
\langle\tilde{\beta}_{\omega},\tilde{\eta}_1\otimes\tilde{\eta}_1\rangle.$$ We now use
$\tilde{\prob},$ the probability measure obtained by transforming
$\prob$ using the functions $\tilde{f}_{\omega}$,
$\tilde{g}_{\omega}$. Note that the colour law in the transformed
measure is now $\eta_1$, and the connectivity
radii~$\tilde{r}_n(a,b)$ satisfy  $$n \, \tilde{r}_n^d(a,b) \to
{\omega(a,b)}/(\eta_1(a)\eta_1(b)) =: \tilde C_d(a,b), \mbox{ as }
n\to\infty.$$ Using (\ref{randomg.Itransform}), we obtain
\begin{align*}
\prob\Big\{(\skril_{\skrig}^1,\skril_{\skrig}^2)& \in
O\Big\}\ge\tilde{\me}\Big\{\sfrac{d\prob}{d\tilde{\prob}}(\skrig)\one_{\{(\skril_{\skrig}^1,\skril_{\skrig}^2)\in
B_{\eta_1}^{1} \times B_{\omega}^{2}\}}\Big\}\\
&=\tilde{\me}\Big\{\prod_{u\in
V}e^{-\tilde{f}_{\omega}(\sigma(x_u))}\prod_{(u,v)\in
E}e^{-\tilde{g}_\omega(\sigma(x_u),\sigma(x_v))}\prod_{(u,v)\in
\skrie}e^{-\sfrac 1n \, \tilde{h}_{\omega, n}(\sigma(x_u),\sigma(x_v))}
\one_{\{(\skril_{\skrig}^1,\skril_{\skrig}^2)\in B_{\eta_1}^{1} \times B_{\omega}^{2}\}}\Big\}\\
&=\tilde{\me}\Big\{ e^{-n\langle \skril_{\skrig}^1,
\tilde{f}_{\omega}\rangle-n\frac{1}{2}\, \langle \skril_{\skrig}^2,
\tilde{g}_{\omega}\rangle-n\frac{1}{2}\,\langle \skril_{\skrig}^1\otimes
\skril_{\skrig}^1, {\tilde{g}}_{\omega}\rangle+\frac{1}{2}\,\langle L_{\Delta}^{1},\tilde{h}_{\omega, n}\rangle}\times\one_{\{(\skril_{\skrig}^1,\skril_{\skrig}^2)\in B_{\eta_1}^{1} \times B_{\omega}^{2}\}}\Big\}\\
&\ge \exp\big(-n\langle \tilde{f}_{\omega},\omega\rangle-n\sfrac 12
\langle\tilde{g}_{\omega},\omega\rangle-n \sfrac 12
\langle\tilde{\beta}_\omega,\eta_1\otimes\eta_1 \rangle
+m-n\eps\big)\times\tilde{\prob}\Big\{(\skril_{\skrig}^1,\skril_{\skrig}^2)\in
B_{\eta_1}^{1} \times B_{\omega}^{2}\Big\},
\end{align*}
where  $m:=0 \wedge \min_{a\in\Sigma}\tilde\beta(a,a).$ Therefore,
by (\ref{randomg.sequence}), we have
\begin{align}\label{randomg.lowbound2}
&\liminf_{n\rightarrow\infty}\sfrac{1}{n}\log\prob\Big\{(\skril_{\skrig}^1,\skril_{\skrig}^2)\in
O\Big\}\nonumber\\
&\ge-\langle \tilde{f}_{\omega},\omega\rangle-\sfrac12 \,
\langle\tilde{g}_{\omega},\omega\rangle-\sfrac 12\,
\langle{\tilde{{\beta}}_{\omega}},
\eta_1\otimes\eta_1 \rangle-\eps+\liminf_{n\rightarrow\infty}\sfrac{1}{n}\log\tilde{\prob}
\Big\{(\skril_{\skrig}^1,\skril_{\skrig}^2)\in B_{\eta_1}^{1} \times
B_{\omega}^{2}\Big\}.\nonumber
\end{align}

The result follows once we prove that
\begin{equation}\label{randomg.translowbound}
\liminf_{n\rightarrow\infty}\sfrac{1}{n}\log\tilde{\prob}\Big\{(\skril_{\skrig}^1,\skril_{\skrig}^2)\in
B_{\eta_1}^{1} \times B_{\omega}^{2}\Big\}= 0.
\end{equation}
We  use the upper bound (but now with the law $\prob$ replaced by
$\tilde{\prob}$) to prove (\ref{randomg.translowbound}). Then we
obtain
\begin{align*}
\limsup_{n\rightarrow\infty}\sfrac{1}{n}\log\tilde\prob\big\{(\skril_{\skrig}^1,\skril_{\skrig}^2)
\in (B_{\eta}^{1} \times B_{\omega}^{2})^{c} \big\}&\le
-\inf_{(\tilde{\rho},\tilde{\omega})\in \tilde{F}}
\tilde{I}(\tilde{\rho},\tilde\omega),
\end{align*}
where $\tilde{F}=(B_{\eta_1}^{1} \times B_{\omega}^{2})^{c}$ and
$\tilde{I}(\tilde{\rho}, \tilde\omega):= H(\tilde\omega \, \| \,
\omega) + \sfrac 12 {\mathfrak{H}_{2}}(\tilde\omega\,\|\,
\tilde{\rho})$. It therefore suffices to show that the infimum is
positive. Suppose for contradiction that there exists a sequence
$(\tilde{\rho}_n,\tilde{\omega}_n)\in\tilde{F}$ with
$\tilde{I}(\tilde{\rho}_n,\tilde\omega_n)\downarrow 0.$ Then, because
$\tilde{I}$ is a good rate function and its level sets are compact,
and by  lower semi-continuity of the mapping
$(\tilde{\rho},\tilde{\omega})\mapsto\tilde{I}(\tilde{\rho},\tilde{\omega})$,
we can construct a limit point $(\tilde{\rho},\tilde\omega)\in\tilde{F}$
with $\tilde{I}(\tilde{\rho},\tilde{\omega})=0$ . By
Lemma~\ref{randomg.Vrate} this implies $H(\tilde{\rho}\,\|\,\eta_1)=0$ and
${\mathfrak{H}_2}(\tilde{\omega}\,\|\,\eta_1)=0$, hence
$\tilde{\rho}=\eta_1,$ and
$\tilde{\omega}=\tilde{C}_d\eta_1\otimes\eta_1=\omega$ contradicting
$(\tilde{\rho},\tilde{\omega})\in\tilde{F}$.
\end{Proof}

\section{Proof of Theorem~\ref{main2}}\label{proofmain2}
For any $n\in\N$ we define
$$\begin{aligned}
\skrip_n(\Sigma) & := \big\{ \rho\in \skrip(\Sigma) \, : \, n\rho(a) \in \N \mbox{ for all } a\in\Sigma\big\},\\
\tilde \skrip_n(\Sigma \times \Sigma) & := \big\{ \omega\in
\tilde\skrip_*(\Sigma\times\Sigma) \, : \, \sfrac
n{1+\one\{a=b\}}\,\omega(a,b) \in \N  \mbox{ for all } a,b\in\Sigma
\big\}\, .
\end{aligned}$$

We denote by
$\Theta_n:=\skrip_n(\Sigma)\times\tilde{\skrip}_n(\Sigma\times\Sigma)$
and
$\Theta:=\skrip(\Sigma)\times\tilde{\skrip}_*(\Sigma\times\Sigma)$.
With
$$\begin{aligned}
P_{(\rho_n, \omega_n)}^{(n)}(\eta_n) & := \prob\big\{\skrim_{\skrig}=\eta_n \, \big| \, \skrih(\skrim_{\skrig})=(\rho_n,\omega_n)\big\}\, ,\\
P^{(n)}(\rho_n,\omega_n) & :=
\prob\big\{(\skril_{\skrig}^1,\skril_{\skrig}^2)=(\rho_n,\omega_n)\big\}
\end{aligned}$$

the joint distribution of $\skril_{\skrig}^1, \skril_{\skrig}^2$ and $\skrim_{\skrig}$ is
the mixture of $P_{(\rho_n, \omega_n)}^{(n)}$ with
$P^{(n)}(\rho_n,\omega_n)$ defined as
\begin{equation}\label{randomg.mixture}
d\tilde{P}^n(\rho_n, \omega_n, \eta_n):= dP_{(\rho_n,
\omega_n)}^{(n)}(\eta_n)\, dP^{(n)}(\rho_n, \omega_n).\,
\end{equation}

(Biggins,2004, Theorem 5(b)) gives criteria for the validity of
large deviation principles for the mixtures and for the goodness of
the rate function if individual large deviation principles are
known. The following three lemmas ensure validity of these
conditions.

We  recall from Lemma~\ref{randomg.uniexpotightness} that the  family of
measures $({P}^n \colon n\in\N)$  is  exponentially tight on
$\Theta$

\begin{lemma}[Doku-Amponsah \& Moerters, 2010] \label{randomg.uniexpotightness} The  family of
measures $(\tilde{P}^n \colon n\in\N)$  is  exponentially tight on
$\Theta\times\skrip(\Sigma\times\N).$
\end{lemma}

Define the function
$$\tilde{J}\colon{\Theta}\times\skrip(\Sigma\times\N)\rightarrow[0,\infty],
\qquad\tilde{J}((\eta_1,\omega),\,\eta)=\tilde{J}_{(\eta_1,\omega)}(\eta),$$
where

\begin{align}\label{randomg.rateLDprob}
\tilde{J}_{(\eta_1,\omega)}(\eta)=\left\{
\begin{array}{ll}H(\eta\,\|\,Q_{poi}^{(\omega, \eta)}) & \mbox
  {if  $(\omega, \eta)$  is  consistent and $\eta_1=\omega_2$ }\\
\infty & \mbox{otherwise.}
\end{array}\right.
\end{align}

\begin{lemma}[Doku-Amponsah \& Moerters, 2010]\label{randomg.convexgoodrate}
$\tilde{J}$ is lower semi-continuous.
\end{lemma}

By (Biggins, 2004, Theorem~5(b)) the two previous lemmas and the
large deviation principles we have established
Theorem~2.2  and (Doku-Amponsah, 2015,  Theorem~2.1) ensure
that under $(\tilde{P}^n)$ the random variables $(\rho_n,
\omega_n, \eta_n)$ satisfy a large deviation principle on
$\skrip(\Sigma) \times \tilde\skrip_*(\Sigma\times\Sigma)
\times\skrip(\Sigma\times\N)$ with good rate function
$$\hat{J}(\eta_1, \omega, \eta)= \left\{ \begin{array}{ll}
H(\eta_1\,\|\,\nu) + \sfrac12\,{\mathfrak H}_2(\omega\,\|\,\Sigma)
+
H(\eta\,\| Q_{poi}^{(\omega,\eta)})\, , & \mbox{ if $(\omega,\eta)$  is  consistent and $\eta_1=\omega_2,$} \\
 \infty\, , & \mbox{ otherwise.} \\ \end{array}\right. $$
By projection onto the last two components we obtain the large
deviation principle as stated in Theorem~\ref{randomg.LDM} from the
contraction principle, see e.g. (Dembo et al.,1998, Theorem~4.2.1).

\section{Proof of Corollary~\ref{ERdd},~Corollary~\ref{RRiv},~and~Corollary~\ref{randomg.L1E}}\label{proofcorollaries}

We derive the theorems from Theorem~\ref{main2} by applying the
contraction principle, see e.g.~(Dembo \& Zeitouni,
1998,  Theorem~4.2.1). In
 fact Theorem~\ref{main2} and the contraction principle imply a large
deviation principle for $D$. It just remains to simplify the rate
functions.

\subsection{Proof of Theorem~\ref{ERdd}.} Note  that,  in the case of an uncoloured
RGG graphs, the function $C$ degenerates to a constant~$c$,
$\skril_{\skrig}^2=|E|/n\in[0,\infty)$ and $\skrim_{\skrig}=D\in\skrip(\N\cup\{0\})$.
Theorem~\ref{main2}  and the contraction principle imply a large
deviation principle for $D$ with good rate function
$$\begin{aligned} \lambda_2(\delta) & =\inf\big\{J(x,\delta) \colon x\ge 0 \big\}
=\inf\big\{ H(\delta\,\|\,q_{x})+\sfrac1 2 x\log x -\sfrac12 x\log
\Delta(d)c+ \sfrac{1}{2}\,\Delta(d)c-\sfrac12 x \colon \langle
\delta\rangle \le x\,\big\},
\end{aligned}$$ which is to be understood as infinity if $\langle
d\rangle$ is infinite. We denote by $\lambda^{x}(\delta)$ the
expression inside the infimum. For any $\eps>0$, we have
$$\begin{aligned}
\lambda_2^{\langle \delta\rangle+\eps}(\delta)-\lambda_2^{\langle
\delta\rangle}(\delta) & =\sfrac{\eps}{2} +\sfrac{\langle
\delta\rangle-\eps}{2}\log\sfrac{\langle \delta\rangle}{\langle
\delta\rangle+\eps}+\sfrac{\eps}{2}\log\sfrac{\langle
\delta\rangle}{\Delta(d)c} \ge \sfrac{\eps}{2} +\sfrac{\langle
\delta\rangle-\eps}{2}\, \big(\sfrac{-\eps}{\langle
\delta\rangle}\big)+\sfrac{\eps}{2}\log\sfrac{\langle
\delta\rangle}{\Delta(d)c}
>0,
\end{aligned}$$
so that the minimum is attained at $x=\Delta(d)\langle\delta\rangle$.

 \subsection{Proof of Corollary~\ref{RRiv}.} Corollary~\ref{ERiv} follows from Theorem~\ref{ERdd} and the
contraction principle applied to the continuous linear map
$G\colon\skrip(\N\cup\{0\})\rightarrow[0,\,1]$ defined by
$G(\delta)=\delta(0).$ Thus, Theorem~\ref{ERdd} implies the large
deviation
 principle for $G(D)=W$ with the good rate function
$\xi_2(y)=\inf\{\lambda_2(\delta) \colon \delta(0)=y, \langle
\delta\rangle < \infty\}.$ We  recall  the  definition  of
$\lambda_2^{x}$ and  observe that $\xi_2(y)$ can be expressed as
$$\xi_2(y)=\inf_{b\ge
0}\inf_{\heap{d\in\skrip(\N\cup\{0\})}{\delta(0)=y,\,
\Delta(d)c\langle \delta\rangle=b^2}} \Big\{\sfrac{1}{2}c+y\log
y+\sfrac{b^2}{2\Delta(d)c} +\sum_{k=1}^{\infty}\delta(k)\log
\sfrac{\delta(k)}{q_{b}(k)}-b(1-y)\Big\}.$$ Now, using Jensen's
inequality, we have that
\begin{equation}\label{randomg.Jensen}\sum_{k=1}^{\infty}\delta(k)\log
\sfrac{\delta(k)}{q_{b}(k)}\ge
(1-y)\log\sfrac{(1-y)}{(1-e^{-b})},\end{equation}
 with equality if $\delta(k)=\sfrac{(1-y)}{(1-e^{-b})}q_{b}(k),$ for all
$k\in\N.$  Therefore, we have the inequality
$$\inf\big\{\lambda_2(\delta) \colon \delta(0)=y,
\langle \delta\rangle <\infty\big\} \ge \inf\big\{
\sfrac{1}{2}c+y\log y+\sfrac{b^2}{2\Delta(d)c}
+(1-y)\log\sfrac{(1-y)}{(1-e^{-b})} -b(1-y) \colon b\ge 0 \big\}.$$
Let $y\in[0,\,1].$  Then, the equation $a(1-e^{-a})=\Delta(d)c(1-y)$
has a unique positive solution. Elementary calculus shows that the
global minimum of
 $b\mapsto\sfrac{1}{2}\Delta(d)c+y\log y+\sfrac{b^2}{2\Delta(d)c} +(1-y)\log\sfrac{(1-y)}{(1-e^{-b})}
-b(1-y)$ on $(0,\infty)$ is attained at the value $b=a$,
 where $a$ is the positive solution of our equation.
We obtain the form of $\xi$  in Corollary~\ref{ERiv} by observing
that
$$\sfrac{a(y)^2+(\Delta(d)c)^2 -2\Delta(d)ca(y)\big(1-y\big)}{2\Delta(d)c}=\sfrac{\Delta(d)cy}{2}\big(2-y\big)+\sfrac{1}{2\Delta(d)c}\big(a(y)-\Delta(d)c(1-y)\big)^2.$$

\subsection{Proof of Corolary~\ref{randomg.L1E}.}
We define the continuous linear map
$W\colon\skrip(\Sigma)\times\tilde{\skrip}_*(\Sigma^2)
\rightarrow[0,\infty)$  by
$W(\eta_1,{\omega})=\sfrac{1}{2}\|{\omega}\|,$ and  infer from
Theorem~\ref{randomg.jointL2L1} and the contraction principle that
$W(\skril_{\skrig}^1,\skril_{\skrig}^2)=|E|/n$ satisfies a large deviation
principle in $[0,\infty)$ with the good rate function
$$\zeta(y)=\inf\big\{I(\eta_1,{\omega}) \colon W(\eta_1,{\omega})=y\big\}.$$
To obtain the form of the rate in the corollary, the infimum  is
reformulated as unconstrained optimization problem (by normalising
$\omega$)
\begin{equation}\label{randomg.Jen}
\inf_{\heap{\omega\in\skrip_*(\Sigma^2)}{\eta_1\in\skrip(\Sigma)}}
\Big\{H(\eta_1\,\|\,\nu)+yH(\omega\,\|\,\Delta(d)C\eta_1\otimes\eta_1 )+y\log2y+
\sfrac{\Delta(d)}{2}\, \| C \omega \otimes \omega\| -y\Big\}.
\end{equation}
By Jensen's inequality
$H(\omega\,\|\,\Delta(d)C\eta_1\otimes\eta_1 )\ge-\log\|\Delta(d)C\eta_1\otimes\eta_1 \|,$
with equality if $\omega=\sfrac{C\eta_1\otimes\eta_1 }{\|
C\eta_1\otimes\eta_1 \|},$ and hence, by symmetry of $C$ we have
\begin{align*}
\min_{\omega\in\skrip_*(\Sigma^2)}& \Big\{
H(\eta_1\,\|\,\nu)+yH(\omega\,\|\,\Delta(d)C\eta_1\otimes\eta_1 )
+y\log2y+\sfrac{\Delta(d)}{2}\, \| C \eta_1\otimes\eta_1 \|-y\Big\}\\
 &= H(\eta_1\,\|\,\nu) -y\log \| \Delta(d)C \eta_1\otimes\eta_1 \|
 +y\log 2y + \sfrac{\Delta(d)}{2}\, \| C \eta_1\otimes\eta_1 \|-y.
\end{align*}

The form given in Corollary~\ref{randomg.L1E} follows by defining
$$y=\sfrac{1}{2}\Delta(d)\sum_{a,b\in\Sigma} C_d(a,b)\eta_1(a)\eta_1(b).$$

\bigskip
{\bf \Large Conclusion}

 In  this  work, we have proved joint large deviation principle for the empirical pair measure
and empirical locality measure of the near intermediate CGRG  models.  From  this  result  we  have  obtained  asymptotic  results  about  useful  graph  quantities  such  as  number of
edges per vertex, the degree distribution and the proportion of isolated vertices for
the near intermediate CGRG  models. The  rate  functions  of  all  these  large  deviation  principles  compared  very well  with  the  rate  functions of the results for  coloured  random  graph  models by  (Doku-Amponsah \&  Moerters,  2010), with  some  extra  terms  accounting  for  the  geometric  effect  in  the  CGRG models. An important future research direction is to formulate and prove an  Asymptotic  Equipartition Property  for Networked Data Structures  Modelled as the  CGRG,  and then  a  possible  Coding  or  Approximate  Pattern Matching Algorithms  for such Networks.  One  could  also  investigate  the   Statistical  Mechanics  on  the  CGRG.\\

{\bf \Large Conflict  of  Interest}

The  author  declares  that  he has  no  conflict  of  interest.\\

{\bf \Large  Acknowledgement}

This  extension  has  been  mentioned  in the author's  PhD Thesis  at  University  of  Bath.


\end{document}